# TIGHT FRAMES AND RELATED GEOMETRIC PROBLEMS

GRIGORY IVANOV

ABSTRACT. We study the properties of a set of vectors called tight frames that obtained as the orthogonal projection of some orthonormal basis of $\mathbb{R}^n$ onto $\mathbb{R}^k$. We show that a set of vectors is a tight frame if and only if the set of all cross products of these vectors is a tight frame. We reformulate a range of problems on the volume of projections (or sections) of regular polytopes in terms of tight frames and write a first-order necessary maximum condition of these problems. As applications, we prove new results for the problem of the maximization of the volume of zonotopes.

## 1. INTRODUCTION

A *tight frame* in $\mathbb{R}^k$ is a set of vectors that is the orthogonal projections of some orthonormal basis of $\mathbb{R}^n$ onto $\mathbb{R}^k$. Or, equivalently, it is a set of vectors $\{v_1, \ldots, v_n\}$ that satisfy the identity

$$\sum_1^n v_i \otimes v_i = I_k, \qquad (1.1)$$

where $I_k$ is the identity operator in $\mathbb{R}^k$.

The tight frames appear and are used naturally in different branches of mathematics: from quantum mechanics and approximation theory (see [1] for details) to classical problems of convex analysis, starting with the well-known John condition [12] for the ellipsoid of maximal volume in a convex body.

We study properties of tight frames from the algebraic point of view. For this purpose we use exterior algebra for studying properties of the exterior powers of the projection operators in Section 3. The necessary definitions on exterior algebra is included in Section 2.1 for completeness. Among other results on projection operators, we prove

**Theorem 1.1.** *A set of vectors $\{v_1, \ldots, v_n\} \subset \mathbb{R}^k$ is a tight frame if and only if all their cross products $[v_{i_1}, \ldots, v_{i_{k-1}}]$ for all $(k-1)$-tuples $\{i_1, \ldots, i_{k-1}\} \in \binom{[n]}{k-1}$ form a tight frame.*

There are standard geometric objects that can be described either in terms of Grassmannian $\operatorname{Gr}(n, k)$ of $k$-dimensional subspace of $\mathbb{R}^n$ or in terms of tight frames. Let us give some examples.

Let $\{v_1, \ldots, v_n\}$ be the orthogonal projection of the standard basis onto $k$-dimensional subspace $H \subset \mathbb{R}^n$. Consider the following classes of convex polytopes:

(1) Projections of the cross-polytope. Let $\Diamond^n = \{x \in \mathbb{R}^n \mid |x_1| + \cdots + |x_n| \leq 1\}$ be the standard cross-polytope. Then the projection of $\Diamond^n$ onto $H$ is the absolute convex hull of the projections of the standard basis of $\mathbb{R}^n$, that is $\operatorname{co}\{\pm v_1, \ldots, \pm v_n\}$.
(2) Sections of cube. Let $\square^n = [-1, 1]^n$ be the standard $n$-dimensional cube. The section $\square^n \cap H$ is the intersections of strips $\{x \in H : |\langle x, v_i \rangle| \leq 1\}$.
(3) Projections of cube. Then the projection of $\square^n$ onto $H$ is the Minkovski sum of segments $[-v_i, v_i]$.

The author was supported by the Swiss National Science Foundation grant 200021_179133.

Identifying $H$ with $\mathbb{R}^k$, we identify the set $\{v_1, \ldots, v_n\}$ with a tight frame. Clearly, the volumes of these polytopes can be considered as functions on the set of tight frames.

Of course, much more types of convex bodies can be described in both ways. We list these three examples since they can be considered as a standard position of special types of centrally symmetric polytopes. For example, any convex centrally symmetric polytope is an affine image of a central section of the cube or a projection of the cross-polytope of some dimension. The problems of finding the minimal and maximal volumes of these types of polytopes are well-studied (see [2], [16], [3], [14]). However, there are still a lot of open questions. For example, the tight upper bound on the volume of a projection of cross-polytope $\Diamond^n$ onto a $k$-dimensional subspace $H$ of $\mathbb{R}^n$ is unknown. The reasonable conjecture is that the upper bound is $\mathrm{vol}_k \Diamond^k$, and it can be considered as a dual to Vaaler's theorem [21] that states that $\mathrm{vol}_k(\square^n \cap H_k) \geq \mathrm{vol}_k \square^k$.

In Section 4 we show how to reformulate problems about extrema of the volume of sections and projections of regular polytopes in the language of tight frames. For example, one can apply our results for listed above types of polytopes. A similar to ours approach was used by Filliman in [8, 7, 6]. However, Filliman preferred to consider the volume functions as the functions on the Grassmannian. Instead of this, we consider functions on the set of tight frames in $\mathbb{R}^k$; extend the domain of such functions to the set of the ordered sets of $n$ vectors of $\mathbb{R}^k$ in a natural way; and give in Lemma 4.3 a necessary and sufficient condition for local extremum of such functions. Thus, we avoid working in the Grassmannians. The only thing we need the Grassmannians for is to prove a first-order approximation formula for a perturbation of a tight frame. As a consequence of Theorem 1.1 and the Cauchy-Binet formula, we obtain.

**Theorem 1.2.** *Let $\{v_1, \ldots, v_n\}$ be a tight frame in $\mathbb{R}^k$, $\{t_i\}_1^n$ be real numbers and $\{x_i\}_1^n$ be an arbitrary set of vectors of $\mathbb{R}^k$. Then*

$$\det \sum_1^n (v_i + t_i x_i) \otimes (v_i + t_i x_i) = 1 + 2 \sum_{i=1}^n t_i \langle x_i, v_i \rangle + o\left(\sqrt{t_1^2 + \ldots + t_n^2}\right).$$

We illustrate our technique on the problem of maximizing function

(1.2) $$F(H) = \mathrm{vol}_k([0,1]^n | H), \quad \dim H = k \leq n.$$

where $[0,1]^n | H$ is the projection of the cube $[0,1]^n$ onto $H$. The projections of the cube $[0,1]^n$ may be considered as a special position of Minkowski sum of finitely many linear segments, so called *zonotopes*.

Using the results of Section 4, we give a complete system of the first-order necessary conditions of a maximizer of (1.2) in Lemma 5.5. It yields the following geometric result.

**Corollary 1.3.** *Let $H$ be a local maximizer of (1.2) and $Q = [0,1]^n | H$. Let $v_i$ denote the orthogonal projection of the standard basis vector $e_i$ onto $H$, for $i \in [n]$, and $Q|v_i^\perp$ be the projection of $Q$ onto the orthogonal complement of line $\{\lambda v_i | \lambda \in \mathbb{R}\}$ in $H$. Then*

$$\mathrm{vol}_{k-1}\left(Q|v_i^\perp\right) = |v_i| \mathrm{vol}_k Q \quad \text{for all} \quad i \in [n].$$

Also we obtain some properties of the maximizers.

**Theorem 1.4.** *Let $H$ be such that the maximum in (1.2) is attained. Let $v_i$ denote the projection of the standard basis vector $e_i$ onto $H$, for $i \in [n]$. Then for any $i, j \in [n]$ the inequality holds*

$$|v_i|^2 \geq (\sqrt{2}-1)|v_j|^2.$$

As a consequence of Theorem 1.4 and McMullen's symmetric formula ([15]), which is

(1.3) $$F(H) = F(H^\perp),$$

where $H^\perp$ is the orthogonal complement of $H$, we prove the following corollary.

**Corollary 1.5.** *Fix $q \in \mathbb{N}$. Let $H_{n-q}$ be a maximizer of (1.2) for $k = n-q, n \geq q$. By $M_n$ and $m_n$ we denote the maximum and the minimum length of the projections of the standard basis's vectors of $\mathbb{R}^n$ onto $H_{n-q}$, respectively. Then*

$$\frac{m_n}{M_n} \to 1 \qquad as \qquad n \to \infty.$$

## 2. Definitions and Preliminaries

We use $C^n$ to denote an $n$-dimensional cube $\{x : 0 \leq x[i] \leq 1\}$ in $\mathbb{R}^n$. Here and throughout the paper $x[i]$ stands for $i$-th coordinate of a vector $x$. As usual, $\{e_i\}_1^n$ is the standard orthonormal basis of $\mathbb{R}^n$. We use $\langle p, x \rangle$ to denote the *value of a linear functional $p$ at a vector $x$*.

For a convex body $K \subset \mathbb{R}^n$ and a $k$-dimensional subspace $H$ of $\mathbb{R}^n$ we denote by $K \cap H$ and $K|H$ the section of $K$ by $H$ and the orthogonal projection of $K$ onto $H$, respectively. For a $k$-dimensional subspace $H$ of $\mathbb{R}^n$ and a convex body $K \subset H$ we denote by $\text{vol}_k K$ the $k$-dimensional volume of $K$. We use $P^H$ to denote the orthogonal projection onto $H$.

For a positive integer $n$, we refer to the set $\{1, 2, \ldots, n\}$ as $[n]$. The set of all $\ell$-element subsets (or $\ell$-*tuple*) of a set $M \subset [n]$ is denoted by $\binom{M}{\ell}$. For two $\ell$-tuples $I \in \binom{[a]}{\ell}, J \in \binom{[b]}{\ell}$, we will use $M_{\{I,J\}}$ to denote the determinant of the corresponding $\ell$ minor of the $a \times b$ matrix $M$. For the sake of convenience, we will write $M_I$ whenever $I = J$. We use $M_{[i,j]}$ to denote the entry in the $i$-th row and the $j$-th column of a matrix $M$.

For the sake of convenience, we denote by $d_S(\{M\})$ the determinant of $\{v_i, i \in M\}$ for an ordered set of vectors $S = \{v_1, \ldots, v_n\} \subset \mathbb{R}^k$ and a set $M$ of $k$ indexes from $[n]$ (the indexes may repeat).

Recall that a *zonotope* in a vector space is the Minkowski sum of finitely many line segments, i.e. $\sum_{i=1}^n [a_i, b_i]$, where $\sum$ stands for the Minkowski sum, and $[a, b]$ means the line segment between points $a$ and $b$.

### 2.1. Exterior algebra.
For the sake of completeness and clarity, all needed definitions from exterior algebra are given here in a brief and non-canonical way. We assume that the equivalence of our definitions to the usual ones is quite obvious. As a proper introduction to multilinear algebra we refer to Greub's book [10].

Let $H$ be a finite dimensional vector space with inner product $\langle \cdot, \cdot \rangle$. We define the vector space $\Lambda^\ell(H)$ as the space of the multilinear skew-symmetric functions on $\ell$ vectors of $H$ with the natural linear structure. The vectors of $\Lambda^\ell(H)$ are called $\ell$-*forms*. As we consider only linear spaces with inner product, we assume that a space $H$ and its dual $H^*$ coincide. This allows us to simplify the following definitions.

For sequence of $\ell$ vectors $\{x_1, \ldots, x_\ell\} \subset H$ we define an $\ell$-form $x_1 \wedge \cdots \wedge x_\ell$ by its evaluation on vectors $\{y_1, \ldots, y_\ell\} \subset H$ given by

$$x_1 \wedge \cdots \wedge x_\ell(y_1, \ldots, y_\ell) = \det M, \qquad \text{where} \qquad M_{[ij]} = \langle x_i, y_j \rangle, \quad i, j \in [\ell].$$

By the properties of the determinant $x_1 \wedge \cdots \wedge x_\ell$ is a multilinear skew-symmetric function on $\ell$ vectors of $H$. For the sake of convenience, given $n$ vectors $(x_i)_1^n$ we denote by $x_L$ the $\ell$-form $x_{i_1} \wedge \cdots \wedge x_{i_\ell}$ for $\ell$-tuple $L = \{i_1, \ldots, i_\ell\} \in \binom{[n]}{\ell}$.

As usual, we use the set of $\ell$-forms $e_{i_1} \wedge \ldots \wedge e_{i_\ell}$, where $\{i_1, \ldots, i_\ell\} \in \binom{[n]}{\ell}$, as the standard basis of $\Lambda^\ell(\mathbb{R}^n)$. We use the lexicographical order on $\ell$-tuples $\{i_1, \ldots, i_\ell\} \in \binom{[n]}{\ell}$ to assign a number to $e_{i_1} \wedge \cdots \wedge e_{i_\ell}$ in this basis.

Recall that for a vector space $H$, an $\ell$-form $w \in \Lambda^\ell(H)$ is said to be *decomposable* if it can be represented in the form $x_1 \wedge \cdots \wedge x_\ell$ for $(x_i)_1^\ell \subset H$. Every $\ell$-form is a linear combination of some decomposable $\ell$-forms (e.g. of the forms from the standard basis), but not all $\ell$-forms are decomposable, e.g. $e_1 \wedge e_2 + e_3 \wedge e_4 \in \Lambda^2(\mathbb{R}^4)$ is not decomposable. A line $l \subset \Lambda^\ell(\mathbb{R}^n)$ has a decomposable directional vector iff $l$ is "generated" by some $\ell$-dimensional subspace $H_\ell \subset \mathbb{R}^n$.

Since the decomposable forms span $\Lambda^\ell(\mathbb{R}^n)$, it is enough to define a linear operator (or the inner product) on $\Lambda^\ell(\mathbb{R}^n)$ only for decomposable $\ell$-forms, by linearity.

Recall that the exterior $\ell$-power of an operator $A$ on $\mathbb{R}^n$ is a linear operator on $\Lambda^\ell(\mathbb{R}^n)$, which is defined on decomposable forms by

$$\wedge^\ell A(x_1 \wedge \cdots \wedge x_\ell) = Ax_1 \wedge \cdots \wedge Ax_\ell.$$

Recall that the inner product on two decomposable $\ell$-forms $a = a_1 \wedge \cdots \wedge a_\ell$ and $b = b_1 \wedge \cdots \wedge b_\ell$ is defined by

$$\langle a, b \rangle = \langle a_1 \wedge \cdots \wedge a_\ell, b_1 \wedge \cdots \wedge b_\ell \rangle = a_1 \wedge \cdots \wedge a_\ell(b_1, \ldots, b_\ell).$$

In this way, once we fixed the inner product on $\mathbb{R}^n$, we fix the inner product on $\Lambda^\ell(\mathbb{R}^n)$. Then for $\ell \in [n]$ we can define a special isometry, the so called *Hodge star operator*, $\star : \Lambda^\ell(\mathbb{R}^n) \to \Lambda^{n-\ell}(\mathbb{R}^n)$ by the equation

$$a \wedge \star(b) = \langle a, b \rangle e_1 \wedge \cdots \wedge e_n,$$

where $a, b$ are $\ell$-forms.

Using this, one can show the following. For a given $a = a_1 \wedge \cdots \wedge a_\ell \in \Lambda^\ell(\mathbb{R}^n)$ and $b = b_1 \wedge \cdots \wedge b_{n-\ell} \in \Lambda^{n-\ell}(\mathbb{R}^n)$, we get

$$(2.1) \qquad \langle a, \star(b) \rangle = (-1)^{\ell(n-\ell)} \det\{a_1, \ldots, a_\ell, b_1, \ldots, b_{n-\ell}\}.$$

For a $k$-dimensional subspace $H_k \subset \mathbb{R}^n$, we use $\wedge^\ell H_k$ to denote the linear hull of forms $x_1 \wedge \cdots \wedge x_\ell$ in $\Lambda^\ell(\mathbb{R}^n)$ such that $\{x_1, \ldots, x_\ell\} \subset H_k$. Since $H_k$ inherits the inner product on $\mathbb{R}^n$, the space $\wedge^\ell H_k$ inherits the inner product of $\Lambda^\ell(\mathbb{R}^n)$, and thus can be identified with the space $\Lambda^\ell(H_k)$ in the tautological way. This allows us to use the Hodge star operator for the spaces $\Lambda^\ell(H_k)$ and $\Lambda^{(k-\ell)}(H_k)$ for $\ell \in [k]$. Also, the identity (2.1) can be rewritten in the following form. For a given $a = a_1 \wedge \cdots \wedge a_\ell \in \Lambda^\ell(H_k)$ and $b = b_1 \wedge \cdots \wedge b_{k-\ell} \in \Lambda^{k-\ell}(H_k)$, we get

$$(2.2) \qquad \langle a, \star(b) \rangle = (-1)^{\ell(k-\ell)} \det\{a_1, \ldots, a_\ell, b_1, \ldots, b_{k-\ell}\},$$

where we understand the determinant as the determinant of $k$-vectors in a $k$-dimensional space $H_k$.

The next Lemma is a straightforward consequence of the definitions. We show that the $\wedge$ product and $\otimes$ product commute for the exterior powers of operators in our case. Usually, definitions of the outer product involve some kind of (anti-)symmetrization, which can give an additional constant when we exchange the $\wedge$ product and $\otimes$ product. For example, we get $-1$ to the proper power for the skew tensor product of graded algebras $\Lambda(E^*)$ and $\Lambda(E)$ (see [10], Section 6.17) in the general case. To avoid misunderstandings we prove the following statement.

**Lemma 2.1.** *Let $A = \sum_{i=1}^{t} v_i \otimes v_i$ for some integer $t$ and vectors $(v_i)_1^t \subset \mathbb{R}^n$. Then $\wedge^\ell A = \sum_{L \in \binom{[t]}{\ell}} v_L \otimes v_L$.*

**Proof.**
We can assume that $\ell \leqslant \dim \mathrm{Lin}\{v_1, \ldots, v_t\} \leqslant t$, otherwise $\wedge^\ell A = 0 = \sum_{L \in \binom{[n]}{\ell}} v_L \otimes v_L$. It suffices to show the identity only for the decomposable $\ell$-forms. Fix $\{x_1, \ldots, x_\ell\} \subset \mathbb{R}^n$. We have

$$\wedge^\ell A(x_1 \wedge \cdots \wedge x_\ell) = Ax_1 \wedge \cdots \wedge Ax_\ell = \sum_{i=1}^{t} \langle v_i, x_1 \rangle v_i \wedge \cdots \wedge \sum_{i=1}^{t} \langle v_i, x_\ell \rangle v_i.$$

By linearity and by skew-symmetry, we expand the last identity and the coefficient at $v_L = v_{i_1} \wedge \cdots \wedge v_{i_\ell}$ is

$$\sum_\sigma \mathrm{sgn}\, \sigma \prod_{j=1}^{\ell} \langle v_{i_j}, x_{\sigma(j)} \rangle,$$

where the summation is taken over all permutations of $[\ell]$. This is the determinant of the matrix $M_{ij} = \langle v_{i_j}, x_i \rangle$, where $i, j \in [\ell]$. By the definition of the inner product, this determinant is just $\langle v_L, (x_1 \wedge \cdots \wedge x_\ell) \rangle$, that is,

$$\wedge^\ell A(x_{i_1} \wedge \cdots \wedge x_{i_\ell}) = \sum_{L \in \binom{[t]}{\ell}} \langle v_L, (x_1 \wedge \cdots \wedge x_\ell) \rangle v_L = \sum_{L \in \binom{[t]}{\ell}} v_L \otimes v_L (x_1 \wedge \cdots \wedge x_\ell).$$

This completes the proof.
□

Identifying $\wedge^{k-1} H_k$ with $\Lambda^{k-1}(H_k)$, we can identify $H_k$ with the $\wedge^{k-1} H_k$ by the Hodge star operator. Then we can define the *cross product* of $k-1$ vectors $\{x_1, \ldots, x_{k-1}\}$ by

$$[x_1, \ldots, x_{k-1}] = \star(x_1 \wedge \cdots \wedge x_{k-1}).$$

Or, in other words, by linearity of the determinant, the cross product $x = [x_1, \ldots, x_{k-1}]$ of $k-1$ vectors $\{x_1, \ldots, x_{k-1}\}$ in a $k$-dimensional space $H_k$ with the fixed inner product $\langle \cdot, \cdot \rangle$ is the vector defined by

$$\langle x, y \rangle = \det(x_1, \ldots, x_{k-1}, y) \quad \text{for all} \quad y \in H_k.$$

## 3. Properties of the exterior power of a projection operator

As a tight frame is a projection of an orthonormal basis, it is natural that its properties are connected with the properties of a projection operator. Moreover, sometimes it is useful to consider a lifting of a tight frame to an orthonormal basis.

In the following trivial Lemma we understand $\mathbb{R}^k \subset \mathbb{R}^n$ as the subspace of vectors, whose last $n - k$ coordinates are zero. For convenience, we will consider $(w_i)_1^n \subset \mathbb{R}^k \subset \mathbb{R}^n$ to be $k$-dimensional vectors.

**Lemma 3.1.** *The following assertions are equivalent:*
  (1) *a set of vectors $\{w_1, \ldots, w_n\} \subset \mathbb{R}^k$ is a tight frame;*
  (2) *there exists an orthonormal basis $\{f_1, \ldots, f_n\}$ of $\mathbb{R}^n$ such that $w_i$ is the orthogonal projection of $f_i$ onto $\mathbb{R}^k$, for any $i \in [n]$;*
  (3) $\mathrm{Lin}\{w_1, \ldots, w_n\} = \mathbb{R}^k$ *and the Gram matrix $\Gamma$ of vectors $\{w_1, \ldots, w_n\} \subset \mathbb{R}^k$ is the matrix of a projection operator from $\mathbb{R}^n$ onto the linear hull of the rows of matrix $M = (w_1, \ldots, w_n)$.*
  (4) *the $k \times n$ matrix $M = (w_1, \ldots, w_n)$ is a sub-matrix of an orthogonal matrix of order $n$.*

The main observation is that $P$ is the Gram matrix of vectors $\{v_1, \ldots, v_n\} \subset H_k$. Since $\langle Pe_i, e_j \rangle = \langle P^2 e_i, e_j \rangle = \langle Pe_i, Pe_j \rangle = \langle v_i, v_j \rangle$, we have that for a fixed $\ell \in [k]$ and $\ell$-tuple $L \subset \binom{[n]}{\ell}$, the corresponding $\ell \times \ell$ sub-matrix of $P$ is the Gram matrix of vectors $(v_i)_{i \in L}$ and the determinant of this Gram matrix is $P_L$. It is well known that for $\ell$ vectors $(w_i)_1^\ell \subset \mathbb{R}^\ell$, their squared determinant is equal to the determinant of their Gram matrix.

**Lemma 3.2.** *Let $P$ be the orthogonal projection from $\mathbb{R}^n$ onto a $k$-dimensional subspace $H_k$. Then for $\wedge^\ell P$, where $1 \leqslant \ell \leqslant k$, we have*
  (1) $\wedge^\ell P$ *is an $\binom{n}{\ell} \times \binom{n}{\ell}$ matrix such that*

(3.1) $$\wedge^\ell P_{[I,J]} = P_{\{I,J\}} \quad \text{for} \quad I, J \in \binom{[n]}{\ell}.$$

  (2) $\wedge^\ell P : \Lambda^\ell(\mathbb{R}^n) \to \Lambda^\ell(\mathbb{R}^n)$ *is the orthogonal projection onto $\wedge^\ell H_k$.*

**Proof.**
1) By definition $\wedge^\ell P$ is an $\binom{n}{\ell} \times \binom{n}{\ell}$ matrix. The identity (3.1) is the consequence of the definition of the scalar product in $\Lambda^\ell(\mathbb{R}^n)$.

2) By the definitions, for any decomposable $\ell$-form $x_1 \wedge \cdots \wedge x_\ell$, we have
$$\wedge^\ell P(x_1 \wedge \cdots \wedge x_\ell) = Px_1 \wedge \cdots \wedge Px_\ell \in \wedge^\ell H_k.$$
By linearity, we have that $\wedge^\ell Px \in \wedge^\ell H_k$ for an arbitrary $\ell$-form $x$. For $x \in \wedge^\ell H_k$, we see that $\wedge^\ell Px = x$. Thus, we showed that $\operatorname{Im} \wedge^\ell P = \wedge^\ell H_k$ and $\left(\wedge^\ell P\right)^2 = \wedge^\ell P$. This completes the proof.
□

In the following Lemma, we understand $\mathbb{R}^k \subset \mathbb{R}^n$ as the subspace of vectors, whose last $n - k$ coordinates are zero. This embedding generates a natural embedding $\Lambda^\ell(\mathbb{R}^k) \subset \Lambda^\ell(\mathbb{R}^n)$ for $\ell \in [n]$.

**Theorem 3.3.** *The following assertions are equivalent for $\{v_1, \ldots, v_n\} \subset \mathbb{R}^k$:*
  *(1) there exists an orthonormal basis $\{f_1, \ldots, f_n\}$ of $\mathbb{R}^n$ such that $v_i$ is the orthogonal projection of $f_i$ onto $\mathbb{R}^k$, for any $i \in [n]$;*
  *(2) $P_k = \sum_1^n v_i \otimes v_i$, where $P_k$ is the projector from $\mathbb{R}^n$ onto $\mathbb{R}^k$*
  *(3) for any fixed $\ell \in [k-1]$ there exists an orthonormal basis $\{f_L\}$ of $\Lambda^\ell(\mathbb{R}^n)$ such that $v_L$ is the orthogonal projection of $f_L$ onto $\Lambda^\ell(\mathbb{R}^k)$, for all $L \in \binom{[n]}{\ell}$.*
  *(4) for any fixed $\ell \in [k-1]$ the following identity is true*
  $$\Lambda^\ell P_k = \sum_{L \in \binom{[n]}{\ell}} v_L \otimes v_L,$$
  *where $\Lambda^\ell P_k$ is the projector from $\Lambda^\ell(\mathbb{R}^n)$ onto $\Lambda^\ell(\mathbb{R}^k)$.*

**Proof.**
By the equivalence (1) ⇔ (2) in Lemma 3.1 and since $I_k$ is the restriction of $P_k$ onto $\mathbb{R}^k$, we have that (1) ⇔ (2).

Identifying $f_i$ and $e_i$ for $i \in [n]$, we identify $\mathbb{R}^k$ with some subspace $H_k \subset \mathbb{R}^n$, $\Lambda^\ell(\mathbb{R}^k) \subset \Lambda^\ell(\mathbb{R}^n)$ with $\wedge^\ell(H^k) \subset \Lambda^\ell(\mathbb{R}^n)$, and $P_k$ with $P$.

Lemma 3.2 says that $\Lambda^\ell P_k$ is exactly $\wedge^\ell P_k$. Thus, identifying $f_i$ and $e_i$ for $i \in [n]$, we identify $\Lambda^\ell P_k$ with $\wedge^\ell P$. Again, since (1) ⇔ (2) in Lemma 3.1, we get that (3) ⇔ (4).

By assertion (2) in Lemma 3.2, we have that (2) ⇒ (4).

Hence, we must show that the implication (2) ⇐ (4) holds to complete the proof.

Let $\{v_1, \ldots, v_n\} \subset H_k$ such that $\wedge^\ell P = \sum_{L \in \binom{[n]}{\ell}} v_L \otimes v_L$ for a fixed $\ell \in [k-1]$. Let us prove that $P = \sum_1^n v_i \otimes v_i$ to complete the proof. Assume the contrary, i.e. $P \neq \sum_1^n v_i \otimes v_i$, and define $A = \sum_1^n v_i \otimes v_i$. By Lemma 2.1, we have that $\wedge^l P = \wedge^l A$. Since the restriction of $A$ on $H_k$ is a positive-semidefinite operator, the restriction of $A$ has an orthonormal basis of eigenvectors in $H_k$. Let us denote these eigenvectors by $x_1, \ldots, x_k$ and the corresponding eigenvalues by $\lambda_1, \ldots, \lambda_k$. Since $\{x_1, \ldots, x_k\} \subset H_k$ and $\wedge^\ell P = \wedge^\ell A$ is the projection onto $\wedge^\ell H_k$, we have that
$$0 \neq x_{i_1} \wedge \cdots \wedge x_{i_\ell} = Px_{i_1} \wedge \cdots \wedge Px_{i_\ell} = \wedge^\ell P(x_{i_1} \wedge \cdots \wedge x_{i_\ell}) =$$
$$\wedge^\ell A(x_{i_1} \wedge \cdots \wedge x_{i_\ell}) = Ax_{i_1} \wedge \cdots \wedge Ax_{i_\ell} = \left(\prod_{j=1}^\ell \lambda_{i_j}\right) \cdot x_{i_1} \wedge \cdots \wedge x_{i_\ell}$$
for any $\ell$-tuple $\{i_1, \ldots, i_\ell\} \subset \binom{[n]}{\ell}$. Therefore, $\left(\prod_{j=1}^\ell \lambda_{i_j}\right) = 1$ for any $\ell$-tuple $\{i_1, \ldots, i_\ell\} \subset \binom{[n]}{\ell}$.
Obviously, since $0 < \ell < k$, it implies that all eigenvalues are one. Hence, the restriction of $A$ onto $H_k$ is the identity in $H_k$. This means that $A = P$.
□

The implication (2) $\Leftarrow$ (4) of Theorem 3.3 still holds for $\ell = k$. By the same arguments as in the proof, we have that (4) implies $\det \left( \sum_1^n v_i \otimes v_i \right) \Big|_{\mathbb{R}^k} = 1$. Thus, we have.

**Corollary 3.4.** *The following assertions are equivalent for $\{v_1, \ldots, v_n\} \subset \mathbb{R}^k \subset \mathbb{R}^n$:*

*(1)* $\det \left( \sum_1^n v_i \otimes v_i \right) \Big|_{\mathbb{R}^k} = 1$.
*(2) the following identity is true*

$$\Lambda^k P_k = \sum_{L \in \binom{[n]}{k}} v_L \otimes v_L.$$

**Proof of Theorem 1.1**
By the definition of the cross product, Theorem 1.1 is the equivalence (1) $\Leftrightarrow$ (3) of Theorem 3.3 with $\ell = k - 1$.
$\square$

By the identity $(\wedge^\ell P)^2 = \wedge^\ell P$, we have

$$P_{\{I,J\}} = \sum_{L \in \binom{[n]}{\ell}} P_{\{I,L\}} P_{\{L,J\}},$$

for a fixed $\ell \in [k]$ and two $\ell$-tuples $I, J \subset \binom{[n]}{\ell}$. But Theorem 3.3 gives us another identity, which connects the squared $\ell$-dimensional volume (that is $P_I$) with the sum of the squared $k$-dimensional volumes.

**Lemma 3.5.** *Let $\{v_1, \ldots, v_n\}$ give a unit decomposition in $H_k$. Then the following identity holds*

$$P_I = \sum_{Q \in \binom{[n]}{k-\ell}; |Q \cap I| = \emptyset} P_{I \cup Q} = \sum_{T \in \binom{[n]}{k} | I \subset T} P_T,$$

*for a fixed $\ell \in [k]$ and an $\ell$-tuple $I \in \binom{[n]}{\ell}$.*

**Proof.**
Now we identify $\wedge^p H_k$ with $\Lambda^p(H_k)$ for a fixed integer $p \in [k]$. Then the Hodge star operator maps a $p$-form $w \in \Lambda^p(H_k)$ to a $(k-p)$-form $\nu \in \Lambda^{k-p}(H_k)$. By Theorem 3.3, we know that the $(k-\ell)$-forms $v_Q$, where $Q \in \binom{[n]}{k-\ell}$, give us a unit decomposition in $\wedge^{(k-\ell)} H_k \equiv \Lambda^{(k-\ell)}(H_k)$. Using the Hodge star operator, we get that the $\ell$-forms $\star(v_Q)$, where $Q \in \binom{[n]}{k-\ell}$, give a unit decomposition in $\wedge^\ell H_k \equiv \Lambda^\ell(H_k)$. By the definition of the inner product of two $\ell$-forms, we have $P_I = \langle v_I, v_I \rangle$. Now we can expand this using properties of the $\ell$-forms $\star(v_Q)$, where $Q \in \binom{[n]}{k-\ell}$. So,

(3.2) $$P_I = \langle v_I, v_I \rangle = \sum_{Q \in \binom{[n]}{k-\ell}} \langle v_I, \star(v_Q) \rangle \langle \star(v_Q), v_I \rangle.$$

By identity (2.2), $\langle v_I, \star(v_Q) \rangle$ is the determinant of $k$ vectors $(v_i)_{i \in I}, (v_q)_{q \in Q}$ of $H_k$. Thus, $\langle v_I, \star(v_Q) \rangle = 0$ whenever $I \cap Q \neq \emptyset$, and $\langle v_I, \star(v_Q) \rangle^2 = P_{I \cup Q}$ whenever $I \cap Q = \emptyset$. Combining these formulas with (3.2), we obtain

$$P_I = \sum_{Q \in \binom{[n]}{k-\ell}; Q \cap I = \emptyset} \langle v_I, \star(v_Q) \rangle \langle \star(v_Q), v_I \rangle = \sum_{Q \in \binom{[n]}{k-\ell}; |Q \cap I| = \emptyset} P_{I \cup Q}.$$

$\square$

3.1. **Lagrange's identity and McMullen's formula.** McMullen's formula (1.3) follows from Shepard's formula [20] for the volume of a $\sum_{i=1}^{n}[0, b_i] \subset \mathbb{R}^k$

$$\text{vol}_k\left(\sum_{i=1}^n [0,b_i]\right) = \sum_{\{i_1,\ldots,i_k\}\in \binom{[n]}{k}} |\det\{b_{i_1},\ldots,b_{i_k}\}|, \qquad (3.3)$$

and identity $|d_S(\{L\})| = |d_S(\{[n]\setminus L\})|$ for a tight frame $S$ and a subset $L$ of $[n]$. By assertion 4 of Lemma 3.1, this identity is equivalent to identity $|U_L| = |U_{[n]\setminus L}|$ for an orthogonal matrix $U$ of rank $n$ and any $L \subset [n]$, which is known as Lagrange's identity. However, it is natural to prove it in terms of properties of a projection operator.

**Lemma 3.6.** *Let $H_{n-k}$ be the orthogonal complement of $H_k$ in $\mathbb{R}^n$. Denote by $P$ and $P_\perp$ the projections onto $H_k$ and $H_{n-k}$, respectively. Then*

$$\star\left(\wedge^k P[x]\right) = \left(\wedge^{(n-k)} P^\perp\right)[\star x]$$

*for any $x \in \Lambda^k(\mathbb{R}^n)$.*

**Proof.**
We use $\pi$ to denote a unit directional vector of the line $\wedge^k H_k$ generated by $H_k$ in $\Lambda^k(\mathbb{R}^n)$. Clearly, $\star\pi$ is the unit directional vector (with a proper orientation) of the line $\wedge^{(n-k)} H_{n-k}$ generated by $H_{n-k}$ in $\Lambda^{n-k}(\mathbb{R}^n)$. For an arbitrary $x \in \Lambda^k(\mathbb{R}^n)$ and $y \in \Lambda^{n-k}(\mathbb{R}^n)$, we have $\wedge^k P[x] = c_1\pi \in \wedge^k H_k$ and $\left(\wedge^{(n-k)} P^\perp\right)[y] = c_2(\star\pi) \in \wedge^{(n-k)} H_{n-k}$ for some scalars $c_1, c_2$. From this and by linearity of operators $\star\left(\wedge^k P\right)$ and $\left(\wedge^{(n-k)} P^\perp\right)\star$, we obtain that $\star\left(\wedge^k P\right)[x] = c\left(\wedge^{(n-k)} P^\perp\right)[\star x]$ for some scalar $c$. From the chain

$$\star\left(\wedge^k P\right)\pi = \star\pi = \left(\wedge^{(n-k)} P^\perp\right)[\star\pi],$$

we conclude that $c = 1$. This completes the proof.
□

So Lagrange's identity implies the following form with a simple geometric meaning.

**Corollary 3.7.** *In the notation of Lemma 3.6, we have*

$$P_L = P^\perp_{[n]\setminus L} \qquad (3.4)$$

*for $L \in \binom{[n]}{k}$.*

## 4. THE TIGHT FRAMES AND RELATED GEOMETRIC PROBLEMS

4.1. **Equivalent formulation.** Tight frame $\{P^H e_1, \ldots, P^H e_n\}$, where $H$ is a $k$-dimensional subspace of $\mathbb{R}^n$ and $P^H$ is the orthogonal projection onto $H$, can be considered as the ordered set of $n$ vectors of $\mathbb{R}^k$. That allows us to avoid working in Grassmannian $\text{Gr}(n,k)$. However, there is an ambiguity in correspondence $H \to \{P^H e_1, \ldots, P^H e_n\}$. Any choice of orthonormal basis of $H$ gives its own tight frame in $\mathbb{R}^k$, all of them are isometric but different from each other. Let us formalize our technique and explain how one can deal with this ambiguity.

We use $\Omega(n,k)$ to denote the set of all tight frames in $\mathbb{R}^k$. We endow $\Omega(n,k)$ with a metric

$$\text{dist}(\{v_1,\ldots,v_n\}, \{w_1,\ldots,w_n\}) = \sqrt{\sum_1^n |v_i - w_i|^2}.$$

The orthogonal group $\text{O}(k)$ acts on $\Omega(n,k)$ in a usual way

$$U\{v_1,\ldots,v_n\} = \{Uv_1,\ldots,Uv_n\} \quad \text{for} \quad U \in \text{O}(k).$$

The equivalence classes of this action are called *O-classes*, and the $O$-class of a tight frame $S$ is denoted as $[S]$. We say that a function $\Psi: \Omega(n,k) \to \mathbb{R}$ is *O-invariant* if it is the constant function on each $O$-class.

There is a one-to-one correspondence between $\text{Gr}(n,k)$ and $\frac{\Omega(n,k)}{\text{O}(k)}$:

- a map $\alpha : \mathrm{Gr}(n,k) \to \frac{\Omega(n,k)}{\mathrm{O}(k)}$ is defined by

$$\alpha(H) = \left[P^H e_1, \ldots, P^H e_n\right] \quad \text{for} \quad H \in \mathrm{Gr}(n,k).$$

- a map $\beta : \frac{\Omega(n,k)}{\mathrm{O}(k)} \to \mathrm{Gr}(n,k)$ is defined as follows. Let $\{v_1, \ldots, v_n\} \in \Omega(n,k)$. By Lemma 3.1, the Gram matrix $P$ of $\{v_1, \ldots, v_n\}$ is the matrix of projection onto a $k$-dimensional subspace $H$. Then $\beta\left(\left[\{v_1, \ldots, v_n\}\right]\right) = H$. Clearly, an orthogonal transformation doesn't change the Gram matrix. Therefore, $\beta$ is defined correctly.

By assertion 4 of Lemma 3.1, the sets of vectors $\{v_1, \cdots, v_n\}$ and $\{Pe_1, \ldots, Pe_n\}$ are isometric. Thus, $\alpha$ and $\beta$ are the inverse functions of each other.

Consequently, there is a one-to-one correspondence between the space of functions $\mathrm{Gr}(n,k) \to \mathbb{R}$ and the space of $O$-invariant functions $\Omega(n,k) \to \mathbb{R}$. For any function $\Psi : \mathrm{Gr}(n,k) \to \mathbb{R}$, we define its frame-function $\gamma(\Psi) : \Omega(k,n) \to \mathbb{R}$ by

$$\gamma(\Psi)(S) = \Psi(\beta([S])).$$

Clearly, the frame-function is $O$-invariant. The argument above showes that $\Psi$ and $\gamma(\Psi)$ have the same global extrema. Let us show that the local extrema of any continuous function $\mathrm{Gr}(n,k) \to \mathbb{R}$ and its frame-function are the same up to factorization by $\mathrm{O}(k)$. Of course, we need to specify a topology or metric on $\mathrm{Gr}(n,k)$. We consider the natural one as a topology of a homogeneous space $(\mathrm{Gr}(n,k) = \mathrm{O}(n)/(\mathrm{O}(k) \times \mathrm{O}(n-k)))$, which is a metric topology generated by metric

$$\mathrm{Dist}(H, H') = \left\|P^H - P^{H'}\right\|,$$

where $\|\cdot\|$ denotes the operator norm.

**Lemma 4.1.** *Let $\Psi : \mathrm{Gr}(n,k) \to \mathbb{R}$ be a continuous function. Then the following assertions are equivalent:*

*(1) $S \in \Omega(n,k)$ is a local extremum of $\gamma(\Psi)$.*
*(2) $[S] \in \Omega(n,k)$ are local extrema of $\gamma(\Psi)$.*
*(3) $\beta([S])$ is a local extremum of $\Psi$.*

**Proof.**
The equivalence of the first two assertion is trivial as $\gamma(\Psi)$ is $O$-invariant. Let us show that they are equivalent to the third one.

We define a metric on $\frac{\Omega(n,k)}{\mathrm{O}(k)}$ as follows

(4.1) $$\rho(O_1, O_2) = \min\{\mathrm{dist}(S_1, S_2) \mid S_1 \in O_1, S_2 \in O_2\}.$$

Obviously, the minimum is attained. Let us check the triangle inequality for arbitrary $O$-classes $O_1, O_2$ and $O_3$. Since $\mathrm{O}(k)$ acts transitively on each $O$-class, there is $S_2 \in O_2$ such that $\rho(O_1, O_2) = \mathrm{dist}(S_1, S_2)$ for any $S_1 \in O_1$. We assume that the minimum in (4.1) is attained on $S_1$ and $S_2$ for $O_1$ and $O_2$, and on $S_2$ and $S_3$ for $O_2$ and $O_3$. Then

$$\rho(O_1, O_2) + \rho(O_2, O_3) = \mathrm{dist}(S_1, S_2) + \mathrm{dist}(S_2, S_3) \geq \mathrm{dist}(S_1, S_3) \geq \rho(O_1, O_3).$$

Clearly, $[S]$ are the local extrema of $\gamma(\Psi)$ iff $[S]$ is a local extremum for function $\Psi_O : \frac{\Omega(n,k)}{\mathrm{O}(k)} \to \mathbb{R}$ defined by $\Psi_O([S]) = \gamma(\Psi(S))$. To complete the proof it suffices to show that $\mathrm{Gr}(n,k)$ and $\frac{\Omega(n,k)}{\mathrm{O}(k)}$ are homeomorphic. This follows from the observation that the operator norm is equivalent to the Hilbert-Schmidt operator norm $\left(\|P^H\|_{HS} = \sqrt{\sum_1^n |P^H e_i|^2}\right)$, and the latter is equivalent to metric $\rho$ on $\frac{\Omega(n,k)}{\mathrm{O}(k)}$.

□

4.2. **Pertubation of frames.** Our main idea of finding local extrema of an $O$-invariant function $\Psi$ is to transform a given tight frame $S$ to a new one $S'$. However, it is not convenient to restrict ourselves to tight frames. For example, consider function $\Psi_1 : \Omega(n, k) \to \mathbb{R}_+$ given by $\Psi_1(\{v_1, \ldots, v_n\}) = \mathrm{vol}_k \mathrm{co}\{\pm v_1, \ldots, \pm v_n\}$ (that is, the frame-function of the volume of projection of the cross-polytope). It is clear what is happening with the convex hull when we move one or several vectors of the tight frame, but we may destroy the tight frame condition with such a transformation. It is not a problem as we can easily extend the domain.

**Definition 4.2.** We will say that the ordered set $S = \{v_1, \ldots, v_n\}$ of $n$ vectors of $\mathbb{R}^k$ is a frame if vectors of $S$ span $\mathbb{R}^k$. We consider the following objects related to frames:
   (1) Operator $\sum_{i \in [n]} w_i \otimes w_i$. We use $A_S$ to denote this operator and the matrix of this operator in the standard basis.
   (2) Operator $B_S = A_S^{-\frac{1}{2}}$.
   (3) A subspace $H_k^S \subset \mathbb{R}^n$ which is the linear hull of the rows of the $k \times n$ matrix $(w_1, \ldots, w_n)$ and the projection operator $P^S$ from $\mathbb{R}^n$ onto $H_k^S$. We use the same notation $P^S$ for the matrix of this operator in the standard basis.

By definition $LS = \{Lv_1, \ldots, Lv_n\}$ for a frame $S = \{v_1, \ldots, v_n\}$ and a linear transformation $L$. The operator $B_S$ is well-defined as the condition $\mathrm{Lin}\, S = \mathbb{R}^k$ implies that $A_S$ is a strictly positive operator. Clearly, $B_S$ maps any frame $S$ to the tight frame:

$$\sum_{i=1}^n B_S v_i \otimes B_S v_i = B_S \left( \sum_{i=1}^n v_i \otimes v_i \right) B_S^T = B_S A_S B_S = I_k.$$

It is easy to see that the metric $\mathrm{dist}(\cdot, \cdot)$ on $\Omega(n, k)$ extends to the set of all frames and a small enough perturbation of frame $S$ yields a small perturbation of $B_S$.

Typically, $O$-invariant functions can be extended to the set of all frames in a tautological way, for example, the extensions of $O$-invariant functions for the objects listed above might be
   (1) $\Psi_1(\{v_1, \ldots, v_n\}) = \mathrm{vol}_k \mathrm{co}\{\pm v_1, \ldots, \pm v_n\}$;
   (2) $\Psi_2(\{v_1, \ldots, v_n\}) = \mathrm{vol}_k \left( \bigcap_1^n \{x \in \mathbb{R}^k : |\langle x, v_i \rangle| \leq 1\} \right)$;
   (3) $\Psi_3(\{v_1, \ldots, v_n\}) = \mathrm{vol}_k \left( \sum_1^n [-v_i, v_i] \right)$.

We prefer to extend the functions as follows

$$(4.2) \qquad \Psi\left(\tilde{S}\right) = \frac{\Psi\left(B_{\tilde{S}} \tilde{S}\right)}{\det B_{\tilde{S}}} \quad \text{for a frame } \tilde{S}.$$

One can see that $\Psi_1$ and $\Psi_3$ satisfy this condition, and $\overline{\Psi}_2 = 1/\Psi_2$ satisfies it as well. It looks as a natural extension at least for problems related to the volume of projections or sections of bodies, because of the positive homogeneity of the volume $\mathrm{vol}_k L(K) = |\det L| \mathrm{vol}_k K$, where $K$ is a convex set and $L$ is a linear transformation.

In order to obtain properties of extremizers, we consider a composition of two transformations:

$$(4.3) \qquad S \xrightarrow{\mathbf{T}} \tilde{S} \xrightarrow{B_{\tilde{S}}} S',$$

where $\mathbf{T}$ will be chosen in a specific way and $B_{\tilde{S}}$ just maps $\tilde{S} = T(S)$ to a new tight frame $S'$. Finally, it is easy to write a necessary and sufficient condition for local extremum of an $O$-invariant function.

**Lemma 4.3.** *A tight frame $S = \{v_1, \ldots, v_n\}$ is a local maximum (or minimum) of a nonnegative $O$-invariant function $\Psi$ iff the following inequality holds for any frame $\tilde{S}$ of open neighbourhood*

$U(S)$ *in the set of all frames*

(4.4) $$\frac{\Psi(\tilde{S})}{\Psi(S)} \leq \sqrt{\det A_{\tilde{S}}} \quad (or \geq).$$

*In addition, $S$ is a global maximum (or minimum) iff inequality (4.4) holds for all frames $\tilde{S}$.*

**Proof.**
As mentioned above, $B_{\tilde{S}}\tilde{S}$ is a tight frame and, by continuity, it belongs to a small enough neighbourhood of $S$ in $\Omega(n,k)$ for small enough $U(S)$. Using these observations and the definition of $B_S$, we have

$$\frac{\Psi(\tilde{S})}{\Psi(S)} = \frac{1}{\det B_{\tilde{S}}} \frac{\Psi(B_{\tilde{S}}\tilde{S})}{\Psi(S)} \leq \frac{1}{\det B_{\tilde{S}}} = \sqrt{\det A_{\tilde{S}}}.$$

The equivalence for a global extremum is trivial.
□

Choosing a proper simple transformation $\mathbb{T}$, ( e.g. scaling one or several vectors, moving one vector to the origin, mapping one vector to another) we may understand the geometric meaning of the left-hand side of (4.4). On the other hand, the determinant in the right-hand side of (4.4) can be calculated directly. In particular, the first-order approximation of the determinant is obtained in Theorem 1.2.

### 4.3. Calculation of determinants.

**Proof of Theorem 1.2**
First of all, by the Cauchy-Binet formula, we have

(4.5) $$\det\left(\sum_1^n v_i \otimes v_i\right) = \sum_{I \in \binom{[n]}{k}} P_I^S.$$

By linearity and identity $P_I^S = (d_S(\{I\}))^2$, the theorem follows from the following lemma.

**Lemma 4.4.** *Let $S = \{v_1, \ldots, v_n\}$ be a tight frame and $S'$ be a frame obtained from $S$ by substitution $v_i \to v_i + tx$, where $t \in \mathbb{R}, x \in \mathbb{R}^k$. Then*

$$\sqrt{\det A_{S'}} = 1 + t\langle v_i, x\rangle + o(t).$$

*That is, $\sqrt{\det A_{S'}}$ is a differentiable function of $t$ at $t = 0$, and the derivative equals $\langle v_i, x\rangle$.*

**Proof.**
Given a substitution $v_i \to v_i + tx$ we modify the minor $P_I^S$ iff $i \in I$. By the properties of minors of $P^S$ and the definition of the Hodge star operator, we have that $P_I^S = (\langle v_i, \star(v_{I\setminus i})\rangle)^2$ if $i \in I$. After the substitution, we get

$$(\langle v_i + tx, \star(v_{I\setminus i})\rangle)^2 = (\langle v_i, \star(v_{I\setminus i})\rangle)^2 + 2t\langle v_i, \star(v_{I\setminus i})\rangle\langle \star(v_{I\setminus i}), x\rangle + o(t).$$

Hence,

$$\det A_{S'} = \sum_{I \in \binom{[n]\setminus i}{k}} P_I^S + \sum_{I \in \binom{[n]}{k}, i \in I} P_I^S + 2t \sum_{I \in \binom{[n]}{k}, i \in I} \langle v_i, \star(v_{I\setminus i})\rangle\langle \star(v_{I\setminus i}), x\rangle + o(t).$$

If $i \in J$ for $J \in \binom{[n]}{k-1}$, then we have $\langle v_i, \star(v_J)\rangle = 0$. Therefore, we can rewrite the last identity

$$\det A_{S'} = \sum_{I \in \binom{[n]}{k}} P_I^S + 2t \sum_{J \in \binom{[n]}{k-1}} \langle v_i, \star(v_J)\rangle\langle \star(v_J), x\rangle + o(t).$$

By Theorem 1.1, we know that the vectors $(\star(v_J))_{J\in\binom{[n]}{k-1}}$ give a unit decomposition in $H_k$. Therefore,
$$\sum_{J\in\binom{[n]}{k-1}} \langle v_i, \star(v_J)\rangle\langle \star(v_J), x\rangle = \langle v_i, x\rangle.$$
Using this and by (4.5), we obtain
$$\sqrt{\det A_{S'}} = \sqrt{\det A_S + 2t\langle v_i, x\rangle + o(t)} = 1 + t\langle v_i, x\rangle + o(t).$$
□

Theorem 1.2 gives the first-order approximation of the left-hand side in inequality (4.4). It implies the first-order necessary condition for a differentiable $O$-invariant function $\Psi$.

**Corollary 4.5.** *Let $\Psi$ be an $O$-invariant function, $S = \{v_1, \ldots, v_n\}$ be a local extremum of $\Psi$ on $\Omega(n,k)$, and $\Psi$ is differentiable at $S$. Then*
$$\Psi(\tilde{S}) = \Psi(S)\left(1 + 2\sum_{i=1}^n t_i\langle x_i, v_i\rangle\right) + o\left(\sqrt{t_1^2 + \ldots + t_n^2}\right)$$
*for an arbitrary $\tilde{S} = \{v_1 + t_1 x_1, \ldots, v_n + t_n x_n\}$.*

By the dimension argument, Corollary 4.5 gives a complete system of the first-order necessary conditions. As a small enough perturbation of a frame is still a frame, we see that the set of frames looks like $\mathbb{R}^{nk}$ locally, and we can speak about local convexity of $O$-invariant functions. By a basic fact from sub-differential calculus, inequality (4.4) implies the differentiability of a locally convex function $\Psi$ at its local maximum.

**Corollary 4.6.** *Let $\Psi$ be an $O$-invariant function, $S = \{v_1, \ldots, v_n\}$ be a local extremum of $\Psi$ on $\Omega(n,k)$, and $\Psi$ is locally convex at $S$. Then $\Psi$ is a differentiable function at $S$ and it satisfies the identity of Corollary 4.5.*

Actually, the local convexity of $\Psi$ is a rather typical situation when we consider the volume of the projection of a polytope in $\mathbb{R}^n$ onto a $k$-dimensional subspace, as the latter is piecewise linear on the corresponding Grassmannian by Theorem 1 in [8]. Another example is the so called *linear parameter systems* (see [17]).

Here is another observation. One can notice that $\operatorname{tr} A_S = \sum_1^n |v_i|^2$ for a tight frame $S = \{v_1, \ldots, v_n\}$. Hence $\operatorname{tr} A_S = k$ for a tight frame $S$. Let $\Omega'(n,k)$ denote a class of frames $S = \{v_1, \ldots, v_n\}$ such that $\sum_1^n |v_i|^2 = k (= \operatorname{tr} A_S)$. The same arguments as in Lemma 4.3 imply.

**Corollary 4.7.** *Let $\Psi$ be an $O$-invariant function. Then*
$$\max_{S\in\Omega'(n,k)} \Psi(S) = \max_{S\in\Omega(n,k)} \Psi(S).$$

**Proof.**
Since $\Omega(n,k) \subset \Omega'(n,k)$, it is enough to show that for any $S' \in \Omega'(n,k), S' \notin \Omega(n,k)$ there exists a tight $S$ such that $\Psi(S') < \Psi(S)$. We put $S = B_{S'}S'$. Then, by Lemma 4.3, it is enough to prove that $\det A_{S'} < 1$ in our case. Considering $A_{S'}$ in the basis of its eigenvectors (in which it is a diagonal operator) and using the inequality of arithmetic and geometric means, we obtain
$$\det A_{S'} \leq \left(\frac{\operatorname{tr} A_{S'}}{k}\right)^k = 1,$$
where equality is attained iff all eigenvalues of $A_{S'}$ are ones. This means, the equality is attained iff $S' \in \Omega(n,k)$. This completes the proof.
□

*Remark* 4.8. Corollary 4.7 was proven in [6] (see Theorem 1) with the same idea but different notation for the projections of regular polytopes, particularly, for the mentioned functions $\Psi_1$ and $\Psi_3$. However, Corollary 4.7 gives a nice property of function $1/\Psi_2$ as well.

4.4. **Lifting of frames.** Vectors $\{v_1, \ldots, v_n\}$ of a tight frame $S \in \Omega(n,k)$ have a nice isometric embedding in $\wedge^{k-1}(H^S)$. For a fixed $i \in [n]$, all $\binom{n}{k-1}$ values $d_S(\{i, L\})$, $L \in \binom{[n]}{k-1}$ can be considered as a vector in $\Lambda^{k-1}(\mathbb{R}^n)$. By definition let $d_S(i)$ be a vector in $\Lambda^{k-1}(\mathbb{R}^n)$ such that its $L$th coordinate in the standard basis of $\Lambda^{k-1}(\mathbb{R}^n)$ is $d_S(\{i, L\})$.

**Lemma 4.9.** *Let $S = \{v_1, \ldots, v_n\}$ be an arbitrary tight frame. Then the vectors $d_S(i), i \in [n]$ lie in $\wedge^{k-1}(H^S) \subset \Lambda^{k-1}(\mathbb{R}^n)$. Moreover, they form a tight frame in $\wedge^{k-1}(H^S)$ and the following identity holds*

$$\langle d_S(i), d_S(j) \rangle = \langle v_i, v_j \rangle, \tag{4.6}$$

*for all $i, j \in [n]$.*

**Proof.**
Let $a_j, j \in [k]$ be the rows of the $k \times n$ matrix $M^S = (v_1, \ldots, v_n)$. By the definition of $H^S$, we know that $a_j \in H^S_k, j \in [k]$ and that they form an orthonormal system in $H^S$. Hence, the $(k-1)$-forms $(a_{[k]\setminus j})_{j=1}^k$ are an orthonormal basis of $\wedge^{k-1}(H^S)$. Consider the $(k-1)$-forms $b_i = \sum_{j=1}^k (-1)^{j+1} v_i[j] \cdot a_{[k]\setminus j}$, where $i \in [n]$. Then, by the definition of the inner product on $\Lambda^{k-1}(\mathbb{R}^n)$ and the Laplace expansion of the determinant, we have

$$b_i[L] = \langle b_i, e_L \rangle = \sum_{j=1}^k (-1)^{j+1} v_i[j] \langle a_{[k]\setminus j}, e_L \rangle = \sum_{j=1}^k (-1)^{j+1} v_i[j] M^S_{[k]\setminus j, L} = M^S_{\{i,L\},\{i,L\}} = d_S(i)[L],$$

for $L \in \binom{[n]}{k-1}$, $i \notin L$. We get $b_i[L] = d_S(i)[L] = 0$ if $i \in L$. Thus, $d_S(i) = b_i \in \Lambda^{k-1}(H^S)$ and, clearly, identity (4.6) holds.
□

Lemma 3.1 allows us to describe all substitutions $w_i \to w_i'$, $i \in [n]$, which preserve the operator $A_S$ for a frame $S = \{w_1, \ldots, w_n\}$. But we need a more suitable geometric description. So, let $S = \{w_1, \ldots, w_n\}$ be a frame and $\ell \in [k]$, $L = \{i_1, \ldots, i_\ell\} \subset \binom{[n]}{\ell}$. Consider a substitution $w_i \to w_i'$, for $i \in L$, and $w_i' = w_i$, for $i \notin L$, which preserves $A_S$, and denote by $S' = \{w_1', \ldots, w_n'\}$ the new frame.

In this notation, we get.

**Lemma 4.10.** *The substitution preserves $A_S$ (i.e. $A_S = A_{S'}$) iff there exists an orthogonal matrix $U$ of rank $\ell$ such that*

$$\left(w_{i_1}', \ldots, w_{i_\ell}'\right) = \left(w_{i_1}, \ldots, w_{i_\ell}\right) U. \tag{4.7}$$

*Additionally, in case $S$ is a tight frame, let $\{f_1, \ldots, f_n\}$ be any orthonormal basis of $\mathbb{R}^n$ given by the assertion (2) of Lemma 3.1. Then the substitution preserves $A_S$ iff the vectors of $S'$ are the projection of an orthonormal basis $\{f_1', \ldots, f_n'\}$ of $\mathbb{R}^n$, which is obtained from $\{f_1, \ldots, f_n\}$ by an orthogonal transformation of $(f_i)_{i \in I}$ in their linear hull.*

**Proof.**
The identity $A_{S'} = A_S$ holds if and only if

$$\sum_{i \in I} w_i \otimes w_i = \sum_{i \in I} w_i' \otimes w_i'.$$

Writing this in a matrix form, we get another equivalent statement that the Gram matrices of the rows of the matrices $\left(w_{i_1}', \ldots, w_{i_\ell}'\right)$ and $\left(w_{i_1}, \ldots, w_{i_\ell}\right)$ are the same. The latter is equivalent to the existence of an isometry of $\mathbb{R}^\ell$, which maps the rows of the first matrix to the rows of the second. This isometry defines an orthogonal matrix $U$, which satisfies the assumptions of the lemma.

In case of a tight frame, extending $U$ as the orthogonal transformation of $\mathrm{Lin}(f_i)_{i \in I}$, we obtain the proper $\{f'_1, \ldots, f'_n\}$.
□

Also, we can reformulate the second claim of Lemma 4.10 in the following equivalent way: Given a block-matrix $M = \begin{pmatrix} A & B \\ C & D \end{pmatrix}$ such that the rows of $\begin{pmatrix} A & B \end{pmatrix}$ are orthonormal and the columns of $\begin{pmatrix} A \\ C \end{pmatrix}$ are orthonormal, then $D$ can be chosen in such a way that $M$ will be an orthogonal matrix.

## 5. Zonotopes and their volume

### 5.1. Definitions and history.
A *zonotope* is the *Minkowski sum* of several segments in $\mathbb{R}^k$. Every zonotope can be represented (up to an affine transformation) as a projection of a higher-dimensional cube $\mathrm{C}^n = [0,1]^n$. For example, it follows from the definition of an operator $B_S$. To get more information about zonotopes, we refer the reader to Zong's book [22, Chapter 2] and to [19].

Several different bounds for the maximum in (1.2) are known. For example, G. D. Chakerian and P. Filliman [4] prove

$$\mathrm{vol}_k(\mathrm{C}^n \,|\, H) \leq \sqrt{\frac{n!}{(n-k)!k!}} \quad \text{and} \quad \mathrm{vol}_k(\mathrm{C}^n \,|\, H) \leq \frac{\omega_{k-1}^k}{\omega_k^{k-1}} \left(\frac{n}{k}\right)^{k/2},$$

where $w_i$ is the volume of the $i$-dimensional Euclidean unit ball. The right-hand side inequality is asymptotically tight, as was shown in [6], the left-hand side inequality is tight in the cases $k = 1, 2, n-2, n-1$. The tight upper and lower bounds in the limit case for the volume of a zonotope in a specific position (even for so called $L_p$-zonoids) were obtained in [14, Theorem 2].

Among different upper bounds, all maximizers of (1.2) were described in the cases $k = 1, 2, n-2, n-1$ and $k = 3, n = 6$ in [6]. It appears that projections of the standard basis vectors onto $H_k$ have the same length whenever $H$ is a maximizer of (1.2) for all cases mentioned above, which gives a rather reasonable conjecture.

**Conjecture 5.1.** *The maximum volume of a projection of the cube $\mathrm{C}^n$ on a subspace is attained when the projections of all edges of the cube have the same length.*

Summarizing the observations of Section 4, we have that the problem to find and to study maximizers of

$$(5.1) \qquad F(S) = \mathrm{vol}_k\left(\sum_1^n [0, v_i]\right), \qquad \text{where} \qquad \{v_1, \ldots, v_n\} \in \Omega(n,k)$$

is equivalent to that of (1.2). We illustrate the developed technique and write the first-order necessary condition of local maximum for (5.1) together with some geometric consequences.

### 5.2. The first-order necessary condition.
We start with showing that $F(S)$ is a differentiable function at its local maximum.

**Lemma 5.2.** *Let $S = \{v_1, \ldots, v_n\}$ be a local maximizer of (5.1) and $S'$ be a frame obtained from $S$ by substitution $v_i \to v_i + tx$, where $t \in \mathbb{R}, x \in \mathbb{R}^k$. Then $\frac{F(S')}{F(S)}$ is a differentiable function of $t$ at $t = 0$ and the derivative equals $\langle v_i, x \rangle$.*

**Proof.**
Given a substitution $v_i \to v_i + tx$, we change $d_S(\{L\})$ for $L \in \binom{[n]}{k}$ in (3.3) iff $i \in L$. Since the determinant is a linear function of each vector, this means that $|d_{S'}(\{L\})|$ as function of $t$ (even as function of $x' = tx$) is a convex function of $t$ (or even of $x' = tx \in \mathbb{R}^k$). Therefore, $F(S')/F(S)$ as function of $t$ (or $x' = tx$) is a convex function. The result follows from Corollary 4.6.
□

**Corollary 5.3.** *Let $S = \{v_1, \ldots, v_n\}$ be a local maximizer of (5.1). Then the vectors of $S$ are in general position in $\mathbb{R}^k$, i.e. $d_S(\{L\}) \neq 0$ for each $L \in \binom{[n]}{k}$.*

**Proof.**
Assume the contrary, that there is a $k$-tuple $J$ such that $d_S(\{J\}) = 0$. As the rank of the vectors of $S$ is $k$, this implies that there is an $L = \{i_1, \ldots, i_k\} \in \binom{[n]}{k}$ such that the vectors $\{v_{i_1}, \ldots, v_{i_{k-1}}\}$ are linearly independent and the vectors $\{v_{i_1}, \ldots, v_{i_{k-1}}, v_{i_k}\}$ are linearly dependent (i.e. $d_S(\{L\}) = 0$). Taking $x \neq 0$ in the orthogonal complement of $\text{Lin}\{v_{i_1}, \ldots, v_{i_{k-1}}\}$ and obtaining $S'$ by a substitution $v_{i_k} \to v_{i_k} + tx$, we get that $F(S')/F(S)$ as well as the absolute value of $d_{S'}(\{L\})$ is not differentiable at $t = 0$. This contradicts Lemma 5.2.
□

**Lemma 5.4.** *Function $F(\cdot)$ is differentiable at a local maximizer of (5.1).*

**Proof.**
Let $S = \{v_1, \ldots, v_n\}$ be a local mazimizer of (5.1), and $S'$ be a frame obtained from $S$ by substitution $v_i \to v_i + t_i x_i$, $i \in [n]$.

Function $|d_{S'}(\{L\})|$ as function of $(t_1, \ldots, t_n)$ is the absolute value of a polynomial of $(t_1, \ldots, t_n)$. Hence, it is differentiable at the point $(t_1, \ldots, t_n)$ whenever $d_{S'}(\{L\}) \neq 0$. By Corollary 5.3, we have $d_S(\{L\}) \neq 0$ for all $k$-tuples. Therefore function

$$F(S') = \sum_{L \in \binom{[n]}{k}} |d_{S'}(\{L\})|$$

as function of $(t_1, \ldots, t_n)$ is differentiable at the origin.
□

Now we show the geometric meaning of identities of Corollary 4.5 in our case. We define a *sign-function* $\sigma_S(i, L)$ by

$$\sigma_S(i, L) = \begin{cases} 1/F(S), & \text{if} \quad d(\{i, L\}) > 0; \\ -1/F(S), & \text{if} \quad d(\{i, L\}) < 0; \\ 0, & \text{if} \quad d(\{i, L\}) = 0, \end{cases}$$

for a frame $S$, $i \in [n]$ and $L \in \binom{[n]}{k-1}$.

In the same way as with the vectors $d_S(i), i \in [n]$, we identify $\sigma_S(i)$ with a vector in $\Lambda^{k-1}(\mathbb{R}^n)$ such that its $L$th coordinate in the standard basis of $\Lambda^{k-1}(\mathbb{R}^n)$ is $\sigma_S(i, L)$. As a direct consequence of Lemma 5.4, we obtain.

**Lemma 5.5.** *Let $S = \{v_1, \ldots, v_n\}$ be a local maximizer of (5.1). Then*

(5.2) $\quad \wedge^{k-1} P^S(\sigma_S(i)) = d_S(i) \quad$ and $\quad \langle \sigma_S(i), d_S(j) \rangle = \langle v_i, v_j \rangle, \quad i, j \in [n].$

**Proof.**
By Lemma 4.9, we have that the vectors $(d_S(i))_{i=1}^n$ give a unit decomposition in $\wedge^{k-1} H_k^S$. Therefore, by Theorem 3.3, it is enough to show that the right-hand side of (5.2) is true.

Fix $i, j \in [n]$. Let $S'$ be a frame obtained from $S$ by the substitution $v_i \to v_i + tv_j$. By Corollary 5.3, we have that $d_S(\{L\}) \neq 0$ for every $k$-tuple $L$. Therefore, $d_{S'}(\{L\}) \neq 0$ and $d_{S'}(\{L\})$ has the same sign as $d_S(\{L\})$ for all $k$-tuples $L$ and for a small enough $t$. Using the substitution $v_i \to v_i + tv_j$, we only change the determinants of type $d_S(\{i, J\})$, where $i \notin J \subset \binom{[n]}{k-1}$. Thus, by the properties of absolute value (for a small enough $t$), we get

$$\frac{F(S')}{F(S)} = \frac{\sum_{L \in \binom{[n]}{k}} |d_{S'}(\{L\})|}{F(S)} = 1 + t \sum_{J \in \binom{[n]}{k-1}, i \notin J} \sigma_S(i, J) d_S(\{j, J\}).$$

Since $\sigma_S(i, J) = 0$ whenever $i \in J$, we have that the coefficient at $t$ in the previous formula is

$$\sum_{J \in \binom{[n]}{k-1}} \sigma_S(i, J) d_S(\{j, J\}).$$

But this is the inner product $\langle \sigma_S(i), d_S(j) \rangle$ of the $(k-1)$-forms $\sigma_S(i)$ and $d_S(j)$ written in the standard basis of $\Lambda^{k-1}(\mathbb{R}^n)$. By Lemma 5.4 and Corollary 4.5, we have that $\langle \sigma_S(i), d_S(j) \rangle = \langle v_i, v_j \rangle$.
□

### 5.3. Proofs of Corollary 1.3, Theorem 1.4 and Corollary 1.5.

Corollary 1.3 is a direct consequence of Lemma 5.5 with a simple geometric meaning.

**Proof of Corollary 1.3.**
By Corollary 5.3, $|v_i| \neq 0$ for all $i \in [n]$. Let $S' = \{v'_1, \ldots, v'_n\}$ be projections of the vectors $\{v_1, \ldots, v_n\}$ onto $v_i^\perp$. Clearly, $v'_i = 0$. By Shepard's formula (3.3) and properties of the determinant, we have

$$\operatorname{vol}_{k-1}\left(Q|v_i^\perp\right) = \sum_{J \in \binom{[n]}{k-1}} |d_{S'}(\{J\})| = \frac{1}{|v_i|} \sum_{J \in \binom{[n]}{k-1}} |d_S(\{i, J\})|.$$

By Lemma 5.5, the latter is $|v_i| \operatorname{vol}_k Q$.
□

The idea of the proof of Theorem 1.4 is the following. We use Lemma 4.10 to rotate vectors by $\pi/4$ (i.e. $v_i \to \cos(\pi/4) v_i - \sin(\pi/4) v_j$, $v_j \to \sin(\pi/4) v_i + \cos(\pi/4) v_j$), and after some simple calculations along with Corollary 4.5, we get the inequality of Theorem 1.4.

**Proof of Theorem 1.4.**
We prove the theorem for a maximizer $S = \{v_1, \ldots, v_n\}$ of (5.1), and we fix $i$ and $j$ in $[n]$. We assume that $|v_i|^2 < |v_j|^2$, otherwise there is nothing to prove. Using Lemma 4.10, we have that the substitution $v_i \to \cos(\pi/4) v_i - \sin(\pi/4) v_j$, $v_j \to \sin(\pi/4) v_i + \cos(\pi/4) v_j$ preserves $A_S$ and the absolute value of $d_S(\{L\})$ for all $L \subset \binom{[n]}{k}$ such that $i, j \in L$. Let $S'$ be the tight frame obtained by this substitution. From the choice of $S$, we have $F(S') \leq F(S)$. Expanding these volumes by formula (3.3) and reducing common determinants, we get

$$\frac{\sqrt{2}}{2} \sum_{J \in \binom{[n]}{k-1}, i,j \notin J} (|d_S(\{i, J\}) - d_S(\{j, J\})| + |d_S(\{i, J\}) + d_S(\{j, J\})|) \leq$$

$$\sum_{J \in \binom{[n]}{k-1}, i,j \notin J} (|d_S(\{i, J\})| + |d_S(\{j, J\})|).$$

By the identity $|a+b| + |a-b| = 2 \max\{|a|, |b|\}$, we obtain that each summand in the left-hand side is at least $2 \max\{|d_S(\{i, J\})|, |d_S(\{i, J\})|\}$, and, consequently, is at least $2 |d_S(\{j, J\})|$. Hence, we show that

$$(\sqrt{2} - 1) \sum_{J \in \binom{[n]}{k-1}, i,j \notin J} |d_S(\{j, J\})| \leq \sum_{J \in \binom{[n]}{k-1}, i,j \notin J} |d_S(\{i, J\})|.$$

For all $(k-1)$-tuples $L$ such that $i \in L$ and $j \notin L$ there is one-to-one correspondence with the set of all $(k-1)$-tuples $J$ such that $i \notin J$ and $j \in J$ given by $L \to (L \setminus \{i\}) \cup \{j\}$. In this case, $|d_S(\{j, L\})| = |d_S(\{i, (L \setminus \{i\}) \cup \{j\}\})|$. Adding all such the determinants to the last inequality, we obtain

$$\sum_{J \in \binom{[n]}{k-1},\, i \notin J} |d_S(\{i, J\})| \geq (\sqrt{2}-1) \sum_{J \in \binom{[n]}{k-1},\, i,j \notin J} |d_S(\{j, J\})| + \sum_{J \in \binom{[n]}{k-1},\, i \in J,\, j \notin J} |d_S(\{j, J\})| \geq$$

$$(\sqrt{2}-1) \left( \sum_{J \in \binom{[n]}{k-1},\, i,j \notin J} |d_S(\{j, J\})| + \sum_{J \in \binom{[n]}{k-1},\, i \in J,\, j \notin J} |d_S(\{j, J\})| \right) = (\sqrt{2}-1) \sum_{J \in \binom{[n]}{k-1},\, j \notin J} |d_S(\{j, J\})|.$$

Finally, the sum in the left-hand side is exactly $F(S)\langle \sigma_S(i), d_S(i) \rangle$, and by Lemma 5.2, it is equal to $F(S)|v_i|^2$. Similarly, we have $(\sqrt{2}-1)F(S)|v_j|^2$ in the right-hand side of the last inequality. Dividing by $F(S)$, we obtain $|v_i|^2 \geq (\sqrt{2}-1)|v_j|^2$.
□

As mentioned in the Introduction, Corollary 1.5 is a consequence of Theorem 1.4 and McMullen's symmetric formula (1.3).

**Proof of Corollary 1.5.**

In the notation of the corollary, let $H_q$ be the orthogonal complement of $H_{n-q}$ in $\mathbb{R}^n$. Let $v_i$ and $v_i'$ be the projections of the vector $e_i$ onto $H_{n-q}$ and $H_q$, respectively. Clearly, $|v_i|^2 + |v_i'|^2 = 1$. By Theorem 1.4, we conclude that $|v_i'|^2$ is at most $1/(\sqrt{2}-1)$ and at least $(\sqrt{2}-1)$ of the average squared length of the projections of the standard basis, which is $q/n$. Therefore,

$$1 \geq \frac{m_n}{M_n} \geq \frac{1 - \frac{1}{\sqrt{2}-1}\frac{q}{n}}{1 - (\sqrt{2}-1)\frac{q}{n}},$$

which tends to 1 as $n$ tends to infinity.
□

Institute of Discrete Mathematics and Geometry, TU Wien, Wiedner Hauptstrasse 8–10/104, A-1040 Wien, Austria

Department of Higher Mathematics, Moscow Institute of Physics and Technology, Institutskii pereulok 9, Dolgoprudny, Moscow region, 141700, Russia.

*E-mail address*: `grigory.ivanov@unifr.ch`